\documentstyle[amsfonts, amssymb, latexsym, epsfig, epic, eepic ,12pt]{amsart}

\newtheorem{theorem}{Theorem}[section]
\newtheorem{lemma}[theorem]{Lemma}
\newtheorem{proposition}[theorem]{Proposition}
\newtheorem{corollary}[theorem]{Corollary}
\newtheorem{definition}[theorem]{Definition}

\def\C{{\mbox{\rm\kern.24em
\vrule width.03em height1.43ex depth-.052ex \kern-.26em C}}}
\def\QSet{\mbox{\rm\kern.24em
\vrule width.03em height1.48ex depth-.051ex \kern-.26em Q}}
\def\Z{{\bf Z}}
\def\R{{\mbox{\rm I\kern-.22em R}}}
\def\P{{\bf P}}
\def\Q{{\bf Q}}
\def\T{{\bf T}}
\def\D{{\bf D}}

\def\size{{\rm size}}

\def\energy{{\rm energy}}

\def\111{\gamma}

\def\be#1{\begin{equation}\label{#1}}
\def\bas{\begin{align*}}
\def\eas{\end{align*}}
\def\bi{\begin{itemize}}
\def\ei{\end{itemize}}
\newenvironment{proof}{\noindent {\bf Proof:} }{\endprf\par}
\def \endprf{\hfill  {\vrule height6pt width6pt depth0pt}\medskip}
\def\emph#1{{\it #1}}

\title{$L^p$ estimates for the biest I. The Walsh case}

\author{Camil Muscalu}
\address{Department of Mathematics, UCLA, Los Angeles CA 90095-1555}
\email{camil@@math.ucla.edu}

\author{Terence Tao}
\address{Department of Mathematics, UCLA, Los Angeles CA 90095-1555 }
\email{tao@@math.ucla.edu}

\author{Christoph Thiele}
\address{Department of Mathematics, UCLA, Los Angeles CA 90095-1555}
\email{thiele@@math.ucla.edu}

\include{psfig}
\begin{document}

\begin{abstract}  
We prove $L^p$ estimates 
(Theorem \ref{main})
for the Walsh model of the ``biest'', a trilinear multiplier with singular symbol.  The corresponding estimates for the Fourier model will be obtained in the sequel \cite{mtt:fourierbiest} of this paper.
\end{abstract}

\maketitle

\section{introduction}

The bilinear Hilbert transform can be written (modulo minor modifications) as
$$
B(f_1,f_2)(x):=\int_{\xi_1<\xi_2}
\widehat{f}_1(\xi_1)\widehat{f}_2(\xi_2) e^{2\pi ix(\xi_1+\xi_2)}\,d\xi_1 d\xi_2,$$
where $f_1, f_2$ are test functions on $\R$ and the Fourier transform is defined by
$$ \hat f(\xi) := \int_\R e^{-2\pi i x \xi} f(x)\ dx.$$
From the work of Lacey and Thiele \cite{laceyt1}, \cite{laceyt2} we have the following $L^p$ estimates on $B$:

\begin{theorem}\label{bht}\cite{laceyt1}, \cite{laceyt2}
$B$ maps $L^p \times L^q \to L^r$ whenever $1 < p, q \leq \infty$, $1/p + 1/q = 1/r$, and $2/3 < r < \infty$.
\end{theorem}

In this paper and the sequel \cite{mtt:fourierbiest} we shall study a trilinear variant $T$ of the bilinear Hilbert transform\footnote{This operator should not be confused with the \emph{trilinear Hilbert transform}, in which the constraint $\xi_1 < \xi_2 < \xi_3$ is replaced by something of the form $\xi_1 + 2 \xi_2 + 3 \xi_3 > 0$.  This operator is just barely beyond the reach of the known multilinear techniques.}, defined by 
\begin{equation}\label{oper}
T(f_1,f_2,f_3)(x):=\int_{\xi_1<\xi_2<\xi_3}
\widehat{f}_1(\xi_1)\widehat{f}_2(\xi_2)\widehat{f}_3(\xi_3)
e^{2\pi ix(\xi_1+\xi_2+\xi_3)}\,d\xi_1 d\xi_2 d\xi_3.
\end{equation}
The operator $T$ arises naturally from WKB expansions of eigenfunctions of one-dimensional Schr\"{o}dinger operators, following the work of Christ and Kiselev \cite{ck}. We discuss this connection further in Appendix I.  For these applications it is of interest to obtain $L^p$ estimates on $T$, especially in the case when the functions $f_j$ are in $L^2$.

From the identity
$$
f_1(x) f_2(x) f_3(x) = \int
\widehat{f}_1(\xi_1)\widehat{f}_2(\xi_2)\widehat{f}_3(\xi_3)
e^{2\pi ix(\xi_1+\xi_2+\xi_3)}\,d\xi_1 d\xi_2 d\xi_3
$$
we see that $T$ has the same homogeneity as the pointwise 
product operator, and hence we expect estimates of H\"older type, i.e., 
$T$ maps $L^{p_1} \times L^{p_2} \times L^{p_3}$ to $L^{p'_4}$ when $1/p'_4 = 1/p_1 + 1/p_2 + 1/p_3$.

%

It is well-known that the operator $B$ has a slightly simpler Walsh model analogue $B_{walsh}$ defined using the Walsh transform instead of the Fourier transform, which we now pause to define.

\begin{definition}\label{walsh-def}
For $l\geq 0$ we define the $l$-th \emph{Walsh function}
$w_l$ by the following recursive formulas

\begin{eqnarray*}
w_0         & := & \chi_{[0,1)}         \\
w_{2l}      & := & w_l(2x)+w_l(2x-1)    \\
w_{2l+1}    & := & w_l(2x)- w_l(2x-1).  
\end{eqnarray*}
\end{definition}

\begin{definition}\label{tile-def}
A \emph{tile} $P$ is a half open rectangle $I_P\times\omega_P$ of area one, such that
$I_P$ and $\omega_P$ are dyadic intervals. If 
$P=[2^{-k}n, 2^{-k}(n+1))\times [2^k l, 2^k(l+1))$ 
is such a tile, we
define the corresponding \emph{Walsh wave packet} $\phi_P$ by
\[\phi_P(x):= 2^{k/2}w_l(2^k x-n).\]
\end{definition}

For each tile $P$, note that $\phi_P$ is supported on $I_P$ and has an $L^2$ norm equal to 1.  Also, observe that $\phi_P$ and $\phi_{P'}$ are orthogonal whenever $P$ and $P'$ are disjoint.

\begin{definition}\label{quartile-def}
A \emph{quartile} $P$ is an open rectangle $I_P\times\omega_P$ of area four, such that
$I_P$ and $\omega_P$ are dyadic intervals . 
For any quartile 
\[ P = [2^{-k}n, 2^{-k}(n+1))\times [2^{k+2} l, 2^{k+2}(l+1))\]
we define the sub-tiles $P_1, P_2, P_3 \subset P$ by
\begin{eqnarray*}
P_1 & := & [2^{-k}n, 2^{-k}(n+1))\times [2^k 4l, 2^k(4l+1)) \\
P_2 & := & [2^{-k}n, 2^{-k}(n+1))\times [2^k(4l+1), 2^k(4l+2)) \\
P_3 & := & [2^{-k}n, 2^{-k}(n+1))\times [2^k(4l+2), 2^k(4l+3)).
\end{eqnarray*}
\end{definition}

\setlength{\unitlength}{0.8mm}
\begin{picture}(190,50)

\put(50,5){\line(1,0){40}}
\put(50,15){\line(1,0){40}}
\put(50,25){\line(1,0){40}}
\put(50,35){\line(1,0){40}}
\put(50,45){\line(1,0){40}}
\put(50,5){\line(0,1){40}}
\put(90,5){\line(0,1){40}}

\put(67,8){$P_1$}
\put(67,18){$P_2$}
\put(67,28){$P_3$}

\end{picture}

\begin{definition}\label{bht-def}  If $\P$ is a finite collection of quartiles, the Walsh Bilinear Hilbert transform $B_{walsh,\P}$ is defined by the formula
$$ B_{walsh,\P}(f_1,f_2) := \sum_{P \in \P} \frac{1}{|I_P|^{1/2}}
\langle f_1, \phi_{P_1} \rangle \langle f_2, \phi_{P_2} \rangle \phi_{P_3}.$$
\end{definition}

From the point of view of the time-frequency phase plane, the Fourier and Walsh models are very similar\footnote{Indeed, the Walsh model can be viewed as the analogue of the Fourier model with the underlying group $\R$ being replaced by $(\Z_2)^{\Z}$.}, however the Walsh model, being dyadic, has several convenient features (such as the ability to localize perfectly in both time and frequency simultaneously) which allow for a clearer and less technical treatment than the Fourier case.  
The following theorem is well known:

\begin{theorem}\label{bht-walsh} 
$B_{walsh,\P}$ maps $L^p \times L^q \to L^r$ whenever $1 < p, q \leq \infty$, $1/p + 1/q = 1/r$, and $2/3 < r < \infty$.  The bounds are uniform in $\P$.
\end{theorem}

Most cases of this theorem were proved in \cite{thiele'}, see
also \cite{thiele''}.
As a by-product of our framework we shall be able to give a self-contained 
proof of this theorem in Section \ref{bht-walsh-sec} of this paper.

Just as the bilinear Hilbert transform $B$ has a Walsh model $B_{walsh,\P}$, the trilinear operator $T$ also has a Walsh model $T_{walsh,\P,\Q}$. 

\begin{definition}\label{twalsh-def}  If $\P$, $\Q$ are two finite collections of quartiles, we define the operator $T_{walsh,\P,\Q}$ by
$$T_{walsh,\P,\Q}:=T'_{walsh,\P,\Q}+T''_{walsh,\P,\Q}$$
where
\bas
T'_{walsh,\P,\Q}(f_1,f_2,f_3)&:=
\sum_{P\in\P}\frac{1}{|I_P|^{1/2}}
\langle B_{walsh,P_1,\Q}(f_1,f_2),\phi_{P_1}\rangle
\langle f_3,\phi_{P_2}\rangle
\phi_{P_3}\\
T''_{walsh,\P,\Q}(f_1,f_2,f_3)&:=
\sum_{P\in\P}\frac{1}{|I_P|^{1/2}}
\langle f_1,\phi_{P_1}\rangle
\langle B_{walsh,P_2,\Q}(f_2,f_3),\phi_{P_2}\rangle
\phi_{P_3},
\end{align*}
where for every tile $P$, $B_{walsh,P,\Q}$ is defined by
$$
B_{walsh,P,\Q}(f_1,f_2):=
\sum_{Q\in\Q;\, \omega_{Q_3}\subseteq \omega_P}
\frac{1}{|I_Q|^{1/2}}
\langle f_1,\phi_{Q_1}\rangle
\langle f_2,\phi_{Q_2}\rangle
\phi_{Q_3}.
$$
\end{definition}

In Appendix II we will explain why the operator $T_{walsh,\P,\Q}$ is the natural Walsh analogue of the Fourier operator $T$.  The operator $B_{walsh,P,\Q}$ can be thought of as the restriction of the operator $B_{walsh,\Q}$ 
to the frequency interval $\omega_P$.  

The main purpose of this paper is to obtain a large set of $L^p$ estimates for the Walsh model $T_{walsh,\P,\Q}$ of $T$.  The operator $T$ itself is a little more technical to handle, and the treatment will be deferred to the sequel \cite{mtt:fourierbiest} of this paper.  

Let us consider now the $3$-dimensional affine hyperspace
\[
S:=\{(\alpha_1,\alpha_2,\alpha_3,\alpha_4)\in\R^4\,
|\,\alpha_1 + \alpha_2 + \alpha_3 + \alpha_4=1\}.
\]

\begin{figure}[htbp]\centering
\psfig{figure=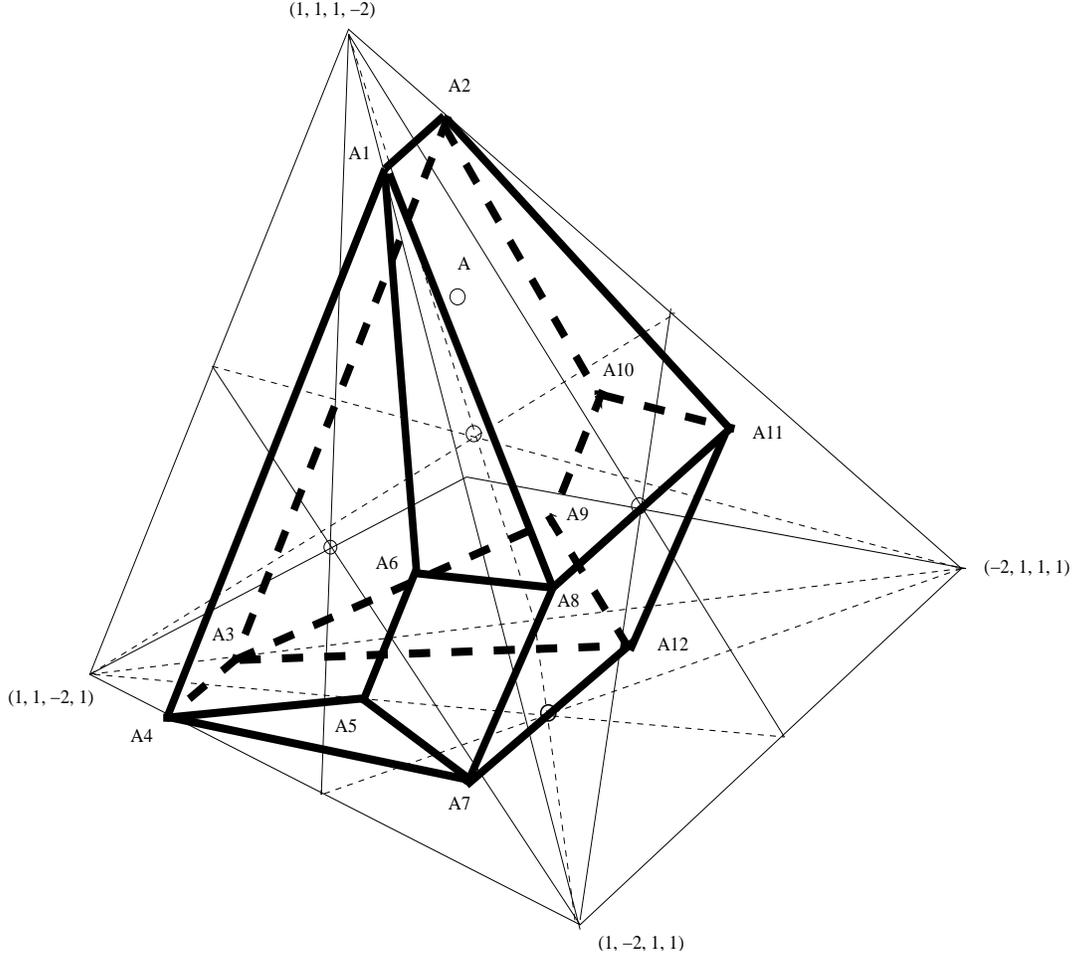, height=5in,width=5.6in}
\caption{The polytope $[A_1...A_{12}]$}
\label{fig}
\end{figure}

Denote by $\D'$ the open interior of the 
convex hull of the $12$ extremal points
$A_1,...,A_{12}$ in Figure 1.
They belong to $S$ and have the following coordinates:

\[
\begin{array}{llll}
A_1:(1,\frac 12,1,-\frac 32)  &  A_2:(\frac 12,1,1,-\frac 32)  &  
A_3:(\frac 12,1,-\frac 32,1) & A_4:(1,\frac 12,-\frac 32,1)  \\  
\ &\ &\ & \ \\
A_5:(1,-\frac 12,0,\frac 12)  &  A_6:(1,-\frac 12,\frac 12,0) & 
A_7:(\frac 12,-\frac 12,0,1)  &  A_8:(\frac 12,-\frac 12,1,0)  \\  
\ &\ &\ & \ \\
A_9:(-\frac 12,1,0,\frac 12) &
A_{10}:(-\frac 12,1,\frac 12,0) &  A_{11}:(-\frac 12,\frac 12,1,0) &  A_{12}:(-\frac 12,\frac 12,0,1).
\end{array}
\]

The point $A$ has the coordinates $(1/2,1/2,1/2,-1/2)$. The other
four circled points, are the centers of gravity of the corresponding
facets of the big tetrahedron, and they have the coordinates
$(1,0,0,0)$, $(0,1,0,0)$, $(0,0,1,0)$, $(0,0,0,1)$. 

Also, denote by $\D''$ the open interior of the 
convex hull of the $12$ extremal points
$B_1,...,B_{12}$ in $S$ (they are not represented in the picture)
where the coordinates of $B_j$ are obtained from the coordinates of
$A_j$ after permuting the indices $1$ and $3$ for $j=1,...,12$ (for
instance $B_2$ has the coordinates $B_2(1,1,1/2,-3/2)$).

Then set $\D:= \D'\cap \D''$.  The (open) 
polytope $\D'$ is the region of estimates
for ${T'}_{walsh,\P,\Q}$, 
while $\D''$ is the region of estimates for 
${T''}_{walsh,\P,\Q}$, as will be clear from the proof of our
main result:

\begin{theorem}\label{main}
Let $1<p_1, p_2, p_3\leq\infty$ and $0<p'_4<\infty$ such that
\[\frac{1}{p_1}+\frac{1}{p_2}+\frac{1}{p_3}=\frac{1}{p'_4}.\]
Then $T_{walsh,\P,\Q}$ maps
\begin{equation}\label{star}
T_{walsh,\P,\Q}: L^{p_1}\times L^{p_2}\times L^{p_3}\rightarrow L^{p'_4}
\end{equation}
with bounds uniform in $\P$, $\Q$
as long as $(1/p_1,1/p_2,1/p_3,1/p_4)\in \D$.
\end{theorem}
In particular,
$T_{walsh,\P,\Q}$ maps $L^{p_1}\times L^{p_2}\times L^{p_3}\rightarrow L^{p'_4}$,
whenever $1<p_1, p_2, p_3\leq\infty$ and $1\leq p'_4 <\infty$.

For the particular application to Schr\"odinger eigenfunctions, the functions $f_i$ will be in $L^2$.  We thus record the corollary (point $A$):
\begin{corollary}\label{23}
$T_{walsh,\P,\Q}$ maps $L^2\times L^2\times L^2\rightarrow L^{2/3}$ uniformly in $\P$, $\Q$.
\end{corollary}

The main new difficulty in treating the operator $T_{walsh,\P,\Q}$, when compared to operators such as $B_{walsh,\P}$, is the presence of the constraint $\omega_{Q_3}\subseteq \omega_{P_i}$ in the double summations.  Without this constraint the operator can be factored into simpler objects.  Eventually we shall exploit the transitivity of the tile ordering $<$ in order to factorize this constraint (see Lemma \ref{biest-trick}).

There are various recurring themes in the subject of multilinear 
singular integrals as in \cite{lacey}, \cite{laceyt1}, \cite{laceyt2}, \cite{cct}
\cite{thiele'}, \cite{thiele}, \cite{gilbertnahmod}, \cite{grafakosli}
and so forth, which the current paper 
again builds up on. While the current article is mostly self contained,
we will mark as ``standard'' any arguments that are well 
understood by now in this framework.

While working on estimating the operator $T$ the authors referred to 
it as the ``beast'' since it seemed worse behaved than previously
seen operators. A teutonic misspelling due to the third author
made it a ``biest'', which then was preferred by the other authors
since it suggests the convenient names ``triest'' or ``multiest'' for the
obvious higher order analogues, which the authors plan to discuss
in forthcoming papers.

The authors would like to thank Mike Christ for pointing out to them 
the occurence of multilinear singular integrals of the type discussed
in this article in the study of eigenfunction expansions of Schr\"odinger operators.

The first author was partially supported by a Sloan Dissertation Fellowship.
The second author is a Clay Prize Fellow and is supported by grants from
the Sloan and Packard Foundations. The third author was partially supported
by a Sloan Fellowship and by NSF grants DMS 9985572 and DMS 9970469.

\section{interpolation}\label{interp-sec}

In this section we review the interpolation theory from \cite{cct} which allows us to reduce multi-linear $L^p$ estimates such as those in Theorem \ref{main} to certain ``restricted weak type'' estimates.

Throughout the paper, we use $A\lesssim B$ to denote the statement that
$A\leq CB$ for some large constant $C$, and $A\ll B$ to denote the
statement that $A\leq C^{-1}B$ for some large constant $C$. Our constants $C$ shall always be independent of $\P$ and $\Q$.

To prove the $L^p$ estimates on $T_{walsh,\P,\Q}$ it is convenient to use duality and introduce the quadrilinear form $\Lambda_{walsh,\P,\Q}$ associated to $T_{walsh,\P,\Q}$ via the formula
$$
\Lambda_{walsh,\P,\Q}(f_1, f_2,f_3, f_4) := \int_{\R} T_{walsh,\P,\Q}(f_1, f_2, f_3)(x)f_4(x) dx.
$$
Similarly define $\Lambda'_{walsh,\P,\Q}$ and $\Lambda''_{walsh,\P,\Q}$.  The statement that $T_{walsh,\P,\Q}$ is bounded from $L^{p_1} \times L^{p_2} \times L^{p_3}$ to $L^{p'_4}$ is then equivalent to $\Lambda_{walsh,\P,\Q}$ being bounded on $L^{p_1} \times L^{p_2} \times L^{p_3} \times L^{p_4}$ if 
$1 <  p_4' < \infty$.  For $p_4'<1$ this simple duality relationship breaks down, however the interpolation arguments in \cite{cct} will allow us to reduce (\ref{star}) to certain ``restricted
type'' estimates on $\Lambda_{walsh,\P,\Q}$. As in \cite{cct} we find more
convenient to work with the quantities $\alpha_i=1/p_i$, $i=1,2,3,4$,
where $p_i$ stands for the exponent of $L^{p_i}$.

\begin{definition}
A tuple $\alpha=(\alpha_1,\alpha_2,\alpha_3,\alpha_4)$ is called admissible, if

\[-\infty <\alpha_i <1\]
for all $1\leq i\leq 4$,

\[\sum_{i=1}^4\alpha_i=1\]
and there is at most one index $j$ such that $\alpha_j<0$. We call
an index $i$ good if $\alpha_i\geq 0$, and we call it bad if
$\alpha_i<0$. A good tuple is an admissible tuple without bad index, a
bad tuple is an admissible tuple with a bad index.
\end{definition}

\begin{definition}
Let $E$, $E'$ be sets of finite measure. We say that $E'$ is a major 
subset of $E$ if $E'\subseteq E$ and $|E'|\geq\frac{1}{2}|E|$.
\end{definition} 

\begin{definition}
If $E$ is a set of finite measure, we denote by $X(E)$ the space of
all functions $f$ supported on $E$ and such that $\|f\|_{\infty}\leq
1$.
\end{definition}

\begin{definition}
If $\alpha=(\alpha_1,\alpha_2,\alpha_3,\alpha_4)$ is an admissible
tuple, we say
that a $4$-linear form $\Lambda$ is of restricted type $\alpha$
if for every sequence $E_1, E_2, E_3, E_4$ of subsets of $\R$ with
finite measure, there exists a major subset $E'_j$ of $E_j$ for each
bad index $j$ (one or none) such that

\[|\Lambda(f_1, f_2, f_3, f_4)|\lesssim |E|^{\alpha}\]
for all functions $f_i\in X(E'_i)$, $i=1,2,3,4$, where we adopt the
convention $E'_i =E_i$ for good indices $i$, and $|E|^{\alpha}$
is a shorthand for

\[|E|^{\alpha}=|E_1|^{\alpha_1}|E_2|^{\alpha_2}
|E_3|^{\alpha_3}|E_4|^{\alpha_4}.\]
\end{definition}

The following ``restricted type'' result will be proved directly.

\begin{theorem}\label{teorema1}
For every vertex $A_i$, $i = 1, \ldots, 12$ there exist
admissible tuples $\alpha$ arbitrarily close to $A_i$ such that the form $\Lambda'_{walsh,\P,\Q}$ is of restricted
type $\alpha$ uniformly in $\P$, $\Q$.
\end{theorem}

By interpolation of restricted weak type estimates (cf. \cite{cct}) we thus obtain

\begin{corollary}\label{restricted}
Let $\alpha$ be an admissible tuple.
Assume also that $\alpha\in \D'$. 
Then $\Lambda'_{walsh,\P,\Q}$ is of restricted type
$\alpha$.
\end{corollary}

Similarly for $\Lambda''$ and $D''$.  Intersecting these two corollaries we obtain the analogous result for $\Lambda$ and $D$.

It only remains to convert these restricted type estimates into strong type estimates.  To do this one just has to apply (exactly as in \cite{cct}) the multilinear Marcinkiewicz interpolation theorem \cite{janson} in the case of good tuples and the interpolation lemma 3.11 in \cite{cct} in the case of bad tuples.

This ends the proof of Theorem \ref{main}. Hence, it remains to prove
Theorem \ref{teorema1}.

\section{trees}

In order to prove the desired estimates for the forms $\Lambda'_{walsh,\P,\Q}$
and $\Lambda''_{walsh,\P,\Q}$ one needs to organize our collections of
quartiles $\P$, $\Q$ into trees as in \cite{fefferman}, 
\cite{laceyt1}, \cite{laceyt2}, \cite{cct}.

\begin{definition}
Let $P$ and $P'$ be tiles. We write $P'<P$ if $I_{P'}\subsetneq I_P$
and $\omega_P\subseteq \omega_{P'}$, and $P'\leq P$ if $P'<P$ or
$P'=P$. 
\end{definition}

Note that $<$ forms a partial order on the set of tiles.  The transitivity of this order shall be crucial, especially in Lemma \ref{biest-trick}.

\begin{definition}
For every $j=1,2,3$ and $P_T\in\P$ define a $j$-tree with top
$P_T$ to be a collection of quartiles $T\subseteq\P$ such that

\[P_j\leq P_{T,j}\]
for all $P\in T$. 
We also say that $T$ is a tree if it is
a $j$-tree for some $j=1,2,3$.
\end{definition}

Notice that $T$ does not necessarily have to contain its top $P_T$.

The following geometric lemma is standard and easy to prove (see 
\cite{laceyt1}, \cite{laceyt2}, \cite{cct}).

\begin{lemma}\label{lacunar}
Let $P, P'$ be quartiles, and $i,j = 1,2,3$ be such that $i \neq j$. If
$P'_i\leq P_i$ then $P'_j\cap P_j=\emptyset$.

In particular, if $T$ is an $i$-tree, then the tiles $\{ P_j: P \in T\}$ are pairwise disjoint.
\end{lemma}

Even more is true: if $T$ is an $i$- tree, then the elements $P$ in
$T$ are parameterized by $I_P$ and the functions $\phi_{P_j}$ behave
like Haar functions in the sense that Calderon- Zygmund theory applies.
Thus an $i$-tree $T$ may be called lacunary in the two indices $j$ other 
than $i$. 

\section{Tile norms}

In the sequel we shall be frequently estimating expressions of the form
\be{trilinear}
| \sum_{P \in \P} \frac{1}{|I_P|^{1/2}} a^{(1)}_{P_1} a^{(2)}_{P_2} a^{(3)}_{P_3}|
\end{equation}
where $\P$ is a collection of quartiles and $a^{(j)}_{P_j}$ are complex numbers for $P \in \P$ and $j = 1,2,3$.  In the treatment of the Walsh bilinear Hilbert transform we just have
\be{aj-def}
a^{(j)}_{P_j} = \langle f_j, \phi_{P_j} \rangle
\end{equation}
 but we will have more sophisticated sequences $a^{(j)}_{P_j}$ when dealing with $\Lambda'_{walsh,\P,\Q}$.  

In order to estimate these expressions it shall be convenient to introduce some norms on sequences of tiles.  The material in this section is standard in the theory of multilinear operators such as the bilinear Hilbert transform, but we reproduce it here for convenience.

\begin{definition}\label{size-def}
Let $\P$ be a finite collection of quartiles, $j=1,2,3$, and let $(a_{P_j})_{P \in \P}$ be a sequence of complex numbers.  We define the \emph{size} of this sequence by
$$ \size_j( (a_{P_j})_{P \in \P} ) := \sup_{T \subset \P}
(\frac{1}{|I_T|} \sum_{P \in T} |a_{P_j}|^2)^{1/2}$$
where $T$ ranges over all trees in $\P$ which are $i$-trees for some $i \neq j$.
We also define the \emph{energy} of the sequence by
$$ \energy_j((a_{P_j})_{P \in \P} ) := \sup_{\D \subset \P}
(\sum_{P \in \D} |a_{P_j}|^2)^{1/2}$$
where $\D$ ranges over all subsets of $\P$ such that the tiles $\{ P_j: P \in \D \}$ are pairwise disjoint.
\end{definition}

The size measures the extent to which the sequence $a_{P_j}$ can concentrate on a single tree and should be thought of as a phase-space variant of the BMO norm.  The energy is a phase-space variant of the $L^2$ norm.  As the  notation suggests, the number $a_{P_j}$ should be thought of as being associated with the tile $P_j$ rather than with the larger quartile $P$.

The usual BMO norm can be written using an $L^2$ oscillation or an $L^1$ oscillation, and the two notions are equivalent thanks to the John-Nirenberg inequality.  The analogous statement for size is

\begin{lemma}\label{jn}  Let $\P$ be a finite collection of quartiles, $j=1,2,3$, and let $(a_{P_j})_{P \in \P}$ be a sequence of complex numbers. Then
\be{jn-est}
\size_j( (a_{P_j})_{P \in \P} ) \sim \sup_{T \subset \P}
\frac{1}{|I_T|} \| ( \sum_{P \in T} |a_{P_j}|^2 \frac{\chi_{I_P}}{|I_P|} )^{1/2} \|_{L^{1,\infty}(I_T)}
\end{equation}
where $T$ ranges over all trees in $\P$ which are $i$-trees for some $i \neq j$.
\end{lemma}

\begin{proof}
Denote the right-hand side of \eqref{jn-est} by $A$.  The bound
$\size_j( (a_{P_j})_{P \in \P} ) \gtrsim A$ is immediate from the H\"older inequality
\bas
\| ( \sum_{P \in T} |a_{P_j}|^2 \frac{\chi_{I_P}}{|I_P|} )^{1/2} \|_{L^{1,\infty}(I_T)}
&\lesssim |I_T|^{1/2}
\| ( \sum_{P \in T} |a_{P_j}|^2 \frac{\chi_{I_P}}{|I_P|} )^{1/2} \|_{L^2(I_T)}\\
&=
|I_T|^{1/2} (\sum_{P \in T} |a_{P_j}|^2)^{1/2}.
\end{align*}

It remains to show $\size_j( (a_{P_j})_{P \in \P} ) \lesssim A$.  We may fix a tree $T$ such that
\be{o-size}
\size_j( (a_{P_j})_{P \in \P} ) = (\frac{1}{|I_T|} \sum_{P \in T} |a_{P_j}|^2)^{1/2}.
\end{equation}
From the definition of $A$ we see that
$$
|\{ ( \sum_{P \in T} |a_{P_j}|^2 \frac{\chi_{I_P}}{|I_P|} )^{1/2} > CA \}|
\leq \frac{1}{2} |I_T|$$
if the constant $C$ is chosen sufficiently large.  

The set on the left-hand side is the union of disjoint dyadic intervals $I_1, \ldots, I_N$.  Observe that
$$
\sum_{P \in T} |a_{P_j}|^2
= (\sum_{i=1}^N \sum_{P \in T: I_P \subseteq I_i} |a_{P_j}|^2)
+ \sum_{P \in T: I_P \not \subseteq I_i \hbox{ for all } 1 \leq i \leq N} |a_{P_j}|^2.$$
By Definition \ref{size-def} we have
$$ \sum_{P \in T: I_P \subseteq I_i} |a_{P_j}|^2 \leq |I_i| \size_j( (a_{P_j})_{P \in \P} )^2.$$
We thus have
$$
\sum_{P \in T} |a_{P_j}|^2
\leq \frac{1}{2} |I_T| | \size_j( (a_{P_j})_{P \in \P} )^2
+ \sum_{P \in T: I_P \not \subseteq I_i \hbox{ for all } 1 \leq i \leq N} |a_{P_j}|^2.$$
From the construction of $I_i$ we have the pointwise estimate
$$ \sum_{P \in T: I_P \not \subseteq I_i \hbox{ for all } 1 \leq i \leq N} |a_{P_j}|^2 \frac{\chi_{I_P}}{I_P} \lesssim A^2.$$
Integrating this on $I_T$ and inserting into the previous we obtain
$$
\sum_{P \in T} |a_{P_j}|^2
\leq \frac{1}{2} |I_T| | \size_j( (a_{P_j})_{P \in \P} )^2 + O(A^2),$$
and the claim follows from (\eqref{o-size}).
\end{proof}

The following estimate is standard but we reproduce a proof in Appendix III for easy reference.  This is the main combinatorial tool needed to obtain estimates on \eqref{trilinear}.

\begin{proposition}\label{abstract}
Let $\P$ be a finite collection of quartiles, and for each $P \in \P$ and $j=1,2,3$ let $a^{(j)}_{P_j}$ be a complex number.  Then
\be{trilinear-est}
| \sum_{P \in \P} \frac{1}{|I_P|^{1/2}} a^{(1)}_{P_1} a^{(2)}_{P_2} a^{(3)}_{P_3}|
\lesssim \prod_{j=1}^3 
\size_j( (a^{(j)}_{P_j})_{P \in \P} )^{\theta_j}
\energy_j( (a^{(j)}_{P_j})_{P \in \P} )^{1-\theta_j}
\end{equation}
for any $0 \leq \theta_1, \theta_2, \theta_3 < 1$ with $\theta_1 + \theta_2 + \theta_3 = 1$, with the implicit constant depending on the $\theta_j$.
\end{proposition}

If we ignore endpoint issues, Proposition \ref{abstract} says that we can estimate \eqref{trilinear} by taking two of the sequences in the energy norm and the third sequence in the size norm.  This is analogous to the H\"older inequality which asserts that a sum $\sum_i a_i b_i c_i$ can be estimated by taking two sequences in $l^2$ and the third in $l^\infty$.

Of course, in order to use Proposition \ref{abstract} we will need some estimates on size and energy.  In the case when $a^{(j)}$ is given by \eqref{aj-def} the relevant estimates are quite straightforward:
 
\begin{lemma}\label{energy-lemma}
Let $j=1,2,3$, $f_j$ be a function in $L^2(\R)$, and let $\P$ be a finite collection of quartiles.  Then we have
\be{energy-lemma-est}
\energy_j((\langle f_j, \phi_{P_j} \rangle)_{P \in \P} ) \leq
\| f_j \|_2.
\end{equation}
\end{lemma}

\begin{proof}  The wave packets $\phi_{P_j}$ are orthonormal whenever the $P_j$ are disjoint.  The claim then follows immediately from Bessel's inequality.
\end{proof}

\begin{lemma}\label{size-lemma}
Let $j=1,2,3$, $E_j$ be a set of finite measure, $f_j$ be a function in $X(E_j)$, and let $\P$ be a finite collection of quartiles.  Then we have
\be{size-lemma-est}
\size_j( (\langle f_j, \phi_{P_j} \rangle)_{P \in \P} ) \lesssim
\sup_{P \in \P} \frac{|E_j \cap I_P|}{|I_P|}.
\end{equation}
\end{lemma}

\begin{proof}
This shall be a Walsh version of the proof of Lemma 7.8 in \cite{cct}.

From Lemma \ref{jn} it suffices to show that
$$
\| F_T \|_{L^{1,\infty}(I_T)} \lesssim  |I_T|  \sup_{P \in \P} \frac{|E_j \cap I_P|}{|I_P|}
$$
for all $i \neq j$ and all trees $T$, where $F_T$ is the vector-valued function
$$ F_T := (\langle f_j, \phi_{P_j} \rangle \frac{\chi_{I_P}}{|I_P|^{1/2}} )_{P \in T}.$$

It suffices to prove this estimate in the case when $T$ contains its top $P_T$, since in the general case one could then decompose $T$ into disjoint trees with this property and then sum.  In this case it thus suffices to show
$$
\| F_T \|_{L^{1,\infty}(I_T)} \lesssim  |E_j \cap I_T|.$$
From the definition of $F_T$ it is clear that we may restrict $f_j$ and $E_j$ to $I_T$, in which case it suffices to show
$$
\| F_T \|_{L^{1,\infty}(I_T)} \lesssim \|f_j\|_1.$$
We shall assume that $T$ is centered at the frequency origin in the sense that 0 is on the boundary of $w_{P_T}$.  (The general case can then be handled by modulating by an appropriate Walsh ``plane wave'').  But then the linear operator $f_j \mapsto F_T$ is a (vector-valued) dyadic Calder\'on-Zygmund operator, and the claim follows from standard theory.
\end{proof}

In the next section we shall show how the above size and energy estimates can be combined with Proposition \ref{abstract} and the interpolation theory of the previous section to obtain Theorem \ref{bht-walsh}.  To prove the estimates for the trilinear operator $T_{walsh,\P,\Q}$ we need some more sophisticated size and energy estimates, which we will pursue after the proof of Theorem \ref{bht-walsh}.

\section{Proof of Theorem \ref{bht-walsh}}\label{bht-walsh-sec}

We now give a proof of Theorem \ref{bht-walsh}.  The proof here is standard, but we give it here for expository purposes, and also because we shall need Theorem \ref{bht-walsh} to prove the size and energy estimates needed for Theorem \ref{teorema1}.

Fix the collection $\P$ of quartiles, and let $\Lambda$ denote the trilinear form
\bas
 \Lambda(f_1,f_2,f_3) &:= \langle B_{walsh,\P}(f_1,f_2),f_3 \rangle\\
&= \sum_{P \in \P} \frac{1}{|I_P|^{1/2}}
\langle f_1, \phi_{P_1} \rangle
\langle f_2, \phi_{P_2} \rangle
\langle f_3, \phi_{P_3} \rangle.
\end{align*}

We shall use the notation of Section \ref{interp-sec}, with the obvious modification for trilinear forms as opposed to quadrilinear forms.
From the interpolation theory in \cite{cct} it suffices to show that $\Lambda$ is of restricted weak type $\alpha$ for all admissible 3-tuples $(\alpha_1, \alpha_2, \alpha_3)$ in the interior of the hexagon with vertices given by the six possible permutations of $(1,1/2,-1/2)$.  By symmetry and interpolation it suffices to prove restricted weak type $\alpha$ for admissible 3-tuples $\alpha$ arbitrarily close to $(1,1/2,-1/2)$, so that the bad index is 3.

Fix $\alpha$ as above, let $E_1$, $E_2$, $E_3$ be sets of finite measure.  We need to find a major subset $E'_3$ of $E_3$ such that
$$
|\Lambda(f_1, f_2, f_3)|\lesssim |E|^{\alpha}
$$
for all functions $f_i\in X(E'_i)$, $i=1,2,3$.

Define the exceptional set $\Omega$ by
\[\Omega := \bigcup_{j=1}^{3}\{M\chi_{E_j}>C |E_j|/|E_3|\}\]
where $M$ is the dyadic Hardy-Littlewood maximal  function.
By the classical Hardy-Littlewood inequality, we have $|\Omega|<1/2 |E_3|$
if $C$ is a sufficiently large constant.  Thus if we set $E'_3 := E_3 \setminus \Omega$, then $E'_3$ is a major subset of $E_3$.

Let $f_i \in X(E'_i)$ for $i=1,2,3$.  We need to show
$$
|\sum_{P \in \P} \frac{1}{|I_P|^{1/2}}
a^{(1)}_{P_1} a^{(2)}_{P_2} a^{(3)}_{P_3}| \lesssim |E|^\alpha$$
where $a^{(j)}_{P_j}$ is defined by \eqref{aj-def}.

We may restrict the quartile set $\P$ to those quartiles $P$ for which $I_P \not \subset \Omega$, since $a^{(3)}_{P_j}$ vanishes for all other quartiles.  By the definition of $\Omega$ we thus have\footnote{Of course, we may also bound the left-hand side trivially by 1.  By combining these two bounds it is possible to prove that $\Lambda$ is of restricted weak type $\alpha$ directly for all $\alpha$ of interest, without recourse to interpolation.}
$$ \frac{|E_j \cap I_P|}{|I_P|} \lesssim \frac{|E_j|}{|E_3|}$$
for all remaining tiles $P \in \P$ and $j=1,2,3$.  From Lemma \ref{size-lemma} we thus have
$$ \size_j( (a^{(j)}_{P_j})_{P \in \P} ) \lesssim \frac{|E_j|}{|E_3|}$$
for $j=1,2,3$.  Also, from Lemma \ref{energy-lemma} and the fact that $f_j \in X(E'_j)$ we have
$$ \energy_j( (a^{(j)}_{P_j})_{P \in \P} ) \lesssim |E_j|^{1/2}.$$
From Proposition \ref{abstract} we thus have
$$
|\sum_{P \in \P} \frac{1}{|I_P|^{1/2}}
a^{(1)}_{P_1} a^{(2)}_{P_2} a^{(3)}_{P_3}| \lesssim 
\prod_{j=1}^3 |E_j|^{(1-\theta_j)/2} (\frac{|E_j|}{|E_3|})^{\theta_j} $$
for any $0 \leq \theta_1, \theta_2, \theta_3 < 1$ such that $\theta_1 + \theta_2 + \theta_3 = 1$.  The claim then follows by choosing $\theta_1 := 2\alpha_1-1$, $\theta_2 := 2\alpha_2-1$, and $\theta_3 := 2\alpha_3 + 1$; note that there exist choices of $\alpha$ arbitrarily close to $(1,1/2,-1/2)$ for which the constraints on $\theta_1, \theta_2, \theta_3, \theta$ are satisfied.  
This concludes the proof of Theorem \ref{bht-walsh}.

\section{Additional size and energy estimates}

We now begin the proof of Theorem \ref{teorema1}.  Fix $\P$, $\Q$
and drop any indices $\P$ and $\Q$ for notational convenience.

In the expression $\Lambda'_{walsh}$ the $Q$ tile in the inner summation has a narrower frequency interval, and hence a wider spatial interval, than the $P$ tile in the outer summation.  Thus the inner summation has a poorer spatial localization than the outer sum.  It shall be convenient to reverse the order of summation so that the inner summation is instead more strongly localized spatially than the outer summation.  Specifically, we rewrite $\Lambda'_{walsh}$ as 
$$
\Lambda'_{walsh}(f_1,f_2,f_3,f_4)=
\sum_{Q\in\Q}\frac{1}{|I_Q|^{1/2}} a^{(1)}_{Q_1} a^{(2)}_{Q_2} a^{(3)}_{Q_3}
$$
where
\begin{equation}\label{rightform}
\begin{split}
a^{(1)}_{Q_1} &:= \langle f_1,\phi_{Q_1} \rangle\\
a^{(2)}_{Q_2} &:= \langle f_2,\phi_{Q_2} \rangle\\
a^{(3)}_{Q_3} &:= \sum_{P\in\P\,;\,\omega_{Q_3}\subseteq \omega_{P_1}}
\frac{1}{|I_P|^{1/2}}
\langle f_3,\phi_{P_2} \rangle
\langle f_4,\phi_{P_3} \rangle
\langle \phi_{P_1}, \phi_{Q_3} \rangle. 
\end{split}
\end{equation}
Observe that for a pairt of quartiles $P$ and $Q$ to give a contribution
to the double sum, we need $P_1\le Q_3$:\\

\setlength{\unitlength}{0.8mm}
\begin{picture}(190,30)

\put(50,10){\line(1,0){40}}
\put(50,10){\line(0,1){10}}
\put(50,20){\line(1,0){40}}
\put(90,10){\line(0,1){10}}

\put(10,15){\line(1,0){160}}
\put(10,15){\line(0,1){2.5}}
\put(10,17.5){\line(1,0){160}}
\put(170,15){\line(0,1){2.5}}

\put(95,8){$P_1$}
\put(175,12){$Q_3$}

\end{picture}

We would like to repeat the argument in Section \ref{bht-walsh-sec}, however we need analogues of Lemma \ref{energy-lemma} and Lemma \ref{size-lemma} for $a^{(3)}_{Q_3}$.  The crucial new ingredient in doing so shall be the following simple geometric lemma which allows us to decouple the $P$ and $Q$ variables.

\begin{lemma}\label{biest-trick}  Let $\D$ be a collection of quartiles such that the tiles $\{ Q_3: Q \in \D \}$ are pairwise disjoint.  Let $\P' \subset \P$ denote the set
$$ \P' := \{ P \in \P: P_1 \leq Q_3 \hbox{ for some } Q \in \D \}.$$
Then for every pair of quartiles $P \in \P$, $Q \in \D$ such that $P_1 \cap Q_3 \neq \emptyset$, we have
$$\omega_{Q_3} \subseteq \omega_{P_1} \hbox{ if and only if } P \in \P'.$$
\end{lemma}

\begin{proof}
Let $P \in \P$, $Q \in \D$ be such that $P_1 \cap Q_3 \neq \emptyset$.

If $\omega_{Q_3} \subseteq \omega_{P_1}$, then $P_1 \leq Q_3$, and so 
$P \in \P'$.  This proves the ``only if'' part.

Now suppose to get a contradiction that there is $P\in \P$ and $R\in \D$
such that $P_1\cap R_3\neq \emptyset$ and 
$\omega_{R_3} \not \subseteq \omega_{P_1}$.  Then $P_1 > R_3$. 
If $P \in \P'$, then we may find $Q \in \D$ such that $P_1 \leq Q_3$, hence 
$R_3 < Q_3$.  But this implies that $R_3 \cap Q_3 \neq \emptyset$, 
contradicting the disjointness hypothesis of the lemma. 
This proves the ``if'' part.
\end{proof}

\setlength{\unitlength}{0.8mm}
\begin{picture}(190,50)

\put(60,10){\line(1,0){10}}
\put(60,10){\line(0,1){40}}
\put(60,50){\line(1,0){10}}
\put(70,10){\line(0,1){40}}

\put(50,30){\line(1,0){40}}
\put(50,30){\line(0,1){10}}
\put(50,40){\line(1,0){40}}
\put(90,30){\line(0,1){10}}

\put(10,35){\line(1,0){160}}
\put(10,35){\line(0,1){2.5}}
\put(10,37.5){\line(1,0){160}}
\put(170,35){\line(0,1){2.5}}

\put(95,28){$P_1$}
\put(75,8){$R_3$}
\put(175,32){$Q_3$}

\multiput(60,35)(2,0){5}{\line(1,1){2}}

\end{picture}

We shall need two analogues of Lemma \ref{energy-lemma}.  The first lemma shall be useful for proving Theorem \ref{teorema1} near the vertices $A_5, \ldots, A_{12}$:

\begin{lemma}\label{bht-energy}
Let $E_j$ be sets of finite measure and $f_j$ be functions in $X(E_j)$ for $j=3,4$.  Then we have
\be{energy-bht-est}
\energy_3((a^{(3)}_{Q_3})_{Q \in \Q}) \lesssim
|E_3|^{(1-\theta)/2} |E_4|^{\theta/2}
\end{equation}
for any $0 < \theta < 1$, with the implicit constant depending on $\theta$.
\end{lemma}

\begin{proof}
By Definition \ref{size-def}, we need to show that
\be{ortho-bht}
(\sum_{Q \in \D} |a^{(3)}_{Q_3}|^2)^{1/2} \lesssim |E_3|^{(1-\theta)/2} |E_4|^{\theta/2}
\end{equation}
for any collection $\D$ of quartiles in $\Q$ such that the tiles $\{ Q_3: Q \in \D\}$ are disjoint.  

Fix $\D$, and define the set $\P'$ by
$$ \P' := \{ P \in \P: P_1 \leq Q_3 \hbox{ for some } Q \in \D \}.$$
By Lemma \ref{biest-trick} and \eqref{rightform} we may write
$$
a^{(3)}_{Q_3} = \sum_{P\in\P'}
\frac{1}{|I_P|^{1/2}}
\langle f_3,\phi_{P_2} \rangle
\langle f_4,\phi_{P_3} \rangle
\langle \phi_{P_1}, \phi_{Q_3} \rangle$$
for all $Q \in \D$.  We can simplify this as
$$ a^{(3)}_{Q_3} = \langle B^*_{walsh,\P'}(f_3,f_4), \phi_{Q_3} \rangle$$
where $B^*_{walsh,\P'}$ is one of the adjoints of the bilinear Hilbert transform $B_{walsh,\P'}$.  Since the $\phi_{Q_3}$ are orthonormal as $Q$ varies in $\D$, we may use Bessel's inequality to estimate the left-hand side of \eqref{ortho-bht} by
$$ \| B^*_{walsh,\P'}(f_3, f_4) \|_2.$$
The claim then follows from Theorem \ref{bht-walsh} and the assumptions $f_3 \in X(E_3)$, $f_4 \in X(E_4)$.
\end{proof}

To prove Theorem \ref{teorema1} near $A_1, \ldots, A_{4}$ we shall use the following sharper variant (The previous lemma follows from this
by the observation $|E_j\cap I_P|\le |I_P|$):

\begin{lemma}\label{bht-energy-2}
Let $E_j$ be sets of finite measure and $f_j$ be functions in $X(E_j)$ for $j=3,4$.  Then we have
\be{energy-bht-est-2}
\energy_3(a^{(3)}_{Q_3}) \lesssim
(|E_4|^{1/2} \sup_{P \in \P} \frac{|E_3 \cap I_P|}{|I_P|})^{1-\theta}
(|E_3|^{1/2} \sup_{P \in \P} \frac{|E_4 \cap I_P|}{|I_P|})^{\theta}
\end{equation}
for any $0 < \theta < 1$, with the implicit constant depending on $\theta$.
\end{lemma}

\begin{proof}
By repeating the proof of Lemma \ref{bht-energy}, we reduce to showing that
$$ \| B^*_{walsh,\P'}(f_3, f_4) \|_2
\lesssim
(|E_4|^{1/2} \sup_{P \in \P} \frac{|E_3 \cap I_P|}{|I_P|})^{1-\theta}
(|E_3|^{1/2} \sup_{P \in \P} \frac{|E_4 \cap I_P|}{|I_P|})^{\theta}
$$
where $\P'$ is an arbitrary subset of $\P$.  By duality we may write the left-hand side as
$$
|\sum_{P \in \P'} \frac{1}{|I_P|^{1/2}}
\langle F, \phi_{P_1} \rangle
\langle f_3, \phi_{P_2} \rangle
\langle f_4, \phi_{P_3} \rangle|
$$
for some $L^2$-normalized function $F$.  By Proposition \ref{abstract} we may estimate this by
\bas
& \energy_1( (\langle F, \phi_{P_1} \rangle)_{P \in \P'})\\
&(\energy_3( (\langle f_4, \phi_{P_3} \rangle)_{P \in \P'})
\size_2( (\langle f_3, \phi_{P_2} \rangle)_{P \in \P'}))^{1-\theta}\\
&(\energy_2( (\langle f_3, \phi_{P_2} \rangle)_{P \in \P'})
\size_3( (\langle f_4, \phi_{P_3} \rangle)_{P \in \P'}))^\theta.
\end{align*}
The claim then follows from Lemma \ref{size-lemma} and Lemma \ref{energy-lemma}.
\end{proof}

The analogue of Lemma \ref{size-lemma} is

\begin{lemma}\label{bht-size}
Let $E_j$ be sets of finite measure and $f_j$ be functions in $X(E_j)$ for $j=3,4$.  Then we have
\be{size-bht-est}
\size_3((a^{(3)}_{Q_3})_{Q \in \Q}) \lesssim
\sup_{Q \in \Q} (\frac{|E_3 \cap I_Q|}{|I_Q|})^{1-\theta} (\frac{|E_4 \cap I_Q|}{|I_Q|})^\theta
\end{equation}
for any $0 < \theta < 1$, with the implicit constant depending on $\theta$.
\end{lemma}

\begin{proof}
By Lemma \ref{jn} it suffices to show that
$$
\| (\sum_{Q \in T} |a^{(3)}_{Q_3}|^2 \frac{\chi_{I_Q}}{|I_Q|})^{1/2} \|_{L^{1,\infty}(I_T)}
\lesssim |I_T|
\sup_{Q \in \Q} (\frac{|E_3 \cap I_Q|}{|I_Q|})^{1-\theta} (\frac{|E_4 \cap I_Q|}{|I_Q|})^\theta
$$
for some $i \neq 3$ and some $i$-tree $T$.  We may assume (as in the proof of Lemma \ref{size-lemma}) that $T$ contains its top $P_T$, in which case we may reduce to
\be{equation}\label{weak}
\| (\sum_{Q \in T} |a^{(3)}_{Q_3}|^2 \frac{\chi_{I_Q}}{|I_Q|})^{1/2} \|_{L^{1,\infty}(I_T)}
\lesssim 
|E_3 \cap I_T|^{1-\theta} |E_4 \cap I_T|^\theta.
\end{equation}
From \eqref{rightform} we see that the only quartiles $P \in \P$ which matter are those such that $I_P \subseteq I_T$.  Thus we may restrict $f_3$, $f_4$, $E_3$, $E_4$ to $I_T$.

Fix $i$, $T$, and define the set $\P'$ by
$$ \P' := \{ P \in \P: P_1 \leq Q_3 \hbox{ for some } Q \in T \}.$$
By Lemma \ref{biest-trick} and \eqref{rightform} as before we have
$$ a^{(3)}_{Q_3} = \langle B^*_{walsh,\P'}(f_3,f_4), \phi_{Q_3} \rangle$$
for all $Q \in T$.  By the dyadic Littlewood-Paley estimate for the tree $T$ we may thus reduce \eqref{weak} to
\[
\| B^*_{walsh,\P'}(f_3,f_4) \|_{L^1(I_T)}
\lesssim 
|E_3|^{1-\theta} |E_4|^\theta.
\]
But this follows from Theorem \ref{bht-walsh} and the assumptions $f_3 \in X(E_3)$, $f_4 \in X(E_4)$.
\end{proof}

\section{Proof of Theorem \ref{teorema1} for $A_5, \ldots, A_{12}$}

Let $\alpha=(\alpha_1,\alpha_2,\alpha_3,\alpha_4)$ admissible tuples near $A_i$ for some $5 \leq i \leq 12$.  We will only consider those vertices with bad index 1 (i.e. $A_9, \ldots, A_{12}$) as the other four vertices can be done similarly.  Thus $\alpha$ has bad index $1$. Let us also
fix $E_1,E_2,E_3,E_4$ arbitrary sets of finite measure. 

As before, we define 
\[\Omega := \bigcup_{j=1}^{4}\{M\chi_{E_j}>C|E_j|/|E_1|\}\]
for a large constant $C$, and set $E'_1 := E_1 \setminus \Omega$.  We now fix $f_i \in X(E'_i)$ for $i=1,2,3,4$.  Our task is then to show

\begin{equation}\label{weakin}
|\sum_{Q\in\Q}\frac{1}{|I_Q|^{1/2}} a^{(1)}_{Q_1} a^{(2)}_{Q_2} a^{(3)}_{Q_3}
|\lesssim |E|^\alpha
\end{equation}
where the $a^{(j)}_{Q_j}$ are defined by \eqref{rightform}.

As before, we may restrict the collection $\Q$ to those quartiles $Q$ for which $I_Q \not \subset \Omega$, since $a^{(1)}_{Q_1}$ vanishes for all other quartiles\footnote{Note however that we cannot restrict $\P$ this way, as 
$I_P\subset I_Q$ and $I_P \subset \Omega$
does not imply $I_Q \subset \Omega$.}.  This implies that
$$ \frac{|E_j \cap I_Q|}{|I_Q|} \lesssim \frac{|E_j|}{|E_1|}$$
for all remaining tiles $Q \in \Q$ and $j=1,2,3,4$.  From Lemma \ref{size-lemma} and Lemma \ref{bht-size} we thus have
\bas
\size_1( (a^{(1)}_{Q_1})_{Q \in \Q} ) &\lesssim \frac{|E_1|}{|E_1|}\\
\size_2( (a^{(2)}_{Q_2})_{Q \in \Q} ) &\lesssim \frac{|E_2|}{|E_1|}\\
\size_3( (a^{(3)}_{Q_3})_{Q \in \Q} ) &\lesssim \frac{|E_3|^{1-\theta} |E_4|^\theta}{|E_1|}
\end{align*}
for some $0 < \theta < 1$ which we will choose later.  Similarly, from Lemma \ref{energy-lemma} and Lemma \ref{bht-energy} and the hypotheses $f_j \in X(E_j)$ we have
\bas
\energy_1( (a^{(1)}_{Q_1})_{Q \in \Q} ) &\lesssim |E_1|^{1/2}\\
\energy_2( (a^{(2)}_{Q_2})_{Q \in \Q} ) &\lesssim |E_2|^{1/2}\\
\energy_3( (a^{(3)}_{Q_3})_{Q \in \Q} ) &\lesssim |E_3|^{(1-\theta)/2} |E_4|^{\theta/2}.
\end{align*}

By Proposition \ref{abstract} we can thus bound the left-hand side of \eqref{weakin} by
$$
\frac{|E_1|^{(1+\theta_1)/2} |E_2|^{(1+\theta_2)/2} (|E_3|^{1-\theta} |E_4|^\theta)^{(1+\theta_3)/2} }{|E_1|},$$
and the claim follows by setting $\theta_1 := 2\alpha_1 + 1$, $\theta_2 := 2\alpha_2 - 1$, $\theta_3 := 2(\alpha_3 + \alpha_4) - 1$, and $\theta := \alpha_4 / (\alpha_3 + \alpha_4)$; the reader may verify that the constraints on $\theta_1, \theta_2, \theta_3, \theta$ can be obeyed for $\alpha$ arbitrarily close to $A_9, A_{10}, A_{11}, A_{12}$.

\section{Proof of Theorem \ref{teorema1} for $A_1, A_2, A_3, A_4$}

Let $\alpha=(\alpha_1,\alpha_2,\alpha_3,\alpha_4)$ admissible tuples near $A_i$ for some $1 \leq i \leq 4$.  We will only consider those vertices with bad index 4 (i.e. $A_1, A_{2}$) as the other two vertices can be done similarly.  Thus $\alpha$ has bad index $4$. Let us also
fix $E_1,E_2,E_3,E_4$ arbitrary sets of finite measure. 

As before, we define 
\[\Omega := \bigcup_{j=1}^{4}\{M\chi_{E_j}>C|E_j|/|E_1|\}\]
for a large constant $C$, and set $E'_4 := E_4 \setminus \Omega$.  We now fix $f_i \in X(E'_i)$ for $i=1,2,3,4$.  Our task is then to show

\begin{equation}\label{weakin-2}
|\sum_{Q\in\Q}\frac{1}{|I_Q|^{1/2}} a^{(1)}_{Q_1} a^{(2)}_{Q_2} a^{(3)}_{Q_3}
|\lesssim |E|^\alpha
\end{equation}
where the $a^{(j)}_{Q_j}$ are defined by \eqref{rightform}.

Recall that $a^{(3)}_{Q_3}$ is defined by
$$
a^{(3)}_{Q_3} := \sum_{P\in\P\,;\,\omega_{Q_3}\subseteq \omega_{P_1}}
\frac{1}{|I_P|^{1/2}}
\langle f_3,\phi_{P_2} \rangle
\langle f_4,\phi_{P_3} \rangle
\langle \phi_{P_1}, \phi_{Q_3} \rangle. 
$$
We may therefore restrict the collection $\P$ to those quartiles $P$ for which $I_P \not \subset \Omega$, since $\langle f_4,\phi_{P_3} \rangle$ vanishes for all other quartiles.  Also observe that $\langle \phi_{P_1}, \phi_{Q_3} \rangle$ vanishes unless $I_P \subset I_Q$.  Thus we may also restrict $\Q$ to those quartiles $Q$ for which $I_Q \not \subset \Omega$.  As a consequence we have 
$$ \frac{|E_j \cap I_Q|}{|I_Q|} \lesssim \frac{|E_j|}{|E_4|}$$
and
$$ \frac{|E_j \cap I_P|}{|I_P|} \lesssim \frac{|E_j|}{|E_4|}$$
for all $P \in \P$, $Q \in \Q$ and $j=1,2,3,4$.  From Lemma \ref{size-lemma} we thus have
\bas
\size_1( (a^{(1)}_{Q_1})_{Q \in \Q} ) &\lesssim \frac{|E_1|}{|E_4|}\\
\size_2( (a^{(2)}_{Q_2})_{Q \in \Q} ) &\lesssim \frac{|E_2|}{|E_4|}.
\end{align*}
From Lemma \ref{bht-size} and the crude estimate $|E_j \cap I_P| \leq |I_P|$ we also have
$$
\size_3( (a^{(3)}_{Q_3})_{Q \in \Q} ) \lesssim 1.$$
Finally, from Lemma \ref{energy-lemma} and Lemma \ref{bht-energy-2} and the hypotheses $f_j \in X(E_j)$ we have
\bas
\energy_1( (a^{(1)}_{Q_1})_{Q \in \Q} ) &\lesssim |E_1|^{1/2}\\
\energy_2( (a^{(2)}_{Q_2})_{Q \in \Q} ) &\lesssim |E_2|^{1/2}\\
\energy_3( (a^{(3)}_{Q_3})_{Q \in \Q} ) &\lesssim 
|E_3|^{(2-\theta)/2} |E_4|^{(\theta-1)/2}
\end{align*}
for some $0 < \theta < 1$ to be chosen later.

By Proposition \ref{abstract} we can thus bound the left-hand side of \eqref{weakin-2} by
$$
\frac{|E_1|^{(1+\theta_1)/2} |E_2|^{(1+\theta_2)/2}}{|E_4|^{1 - \theta_3}} (|E_3|^{(2-\theta)/2} |E_4|^{(\theta-1)/2})^{1-\theta_3},$$
and the claim follows by setting $\theta_1 := 2\alpha_1 - 1$, $\theta_2 := 2\alpha_2 - 1$, $\theta_3 := 2(\alpha_3 + \alpha_4) + 1$, and $\theta := (3\alpha_3 + 2\alpha_4) / (\alpha_3 + \alpha_4)$; the reader may verify that the constraints on $\theta_1, \theta_2, \theta_3, \theta$ can be obeyed for $\alpha$ arbitrarily close to $A_{1}, A_{2}$.

\section{Appendix I: Connection with eigenfunctions of Schr\"odinger operators}

In this section we sketch why the operator $T$ arises naturally from the study of eigenfunctions of Schr\"odinger operators.  Further details can be found in the work of Christ and Kiselev \cite{ck}, \cite{ck-2}.

Let $V(\xi)$ be a locally integrable function on $\R$.  We consider the eigenfunction equation\footnote{In the literature the variable $x$ is usually denoted $k$, while $\xi$ is denoted $x$.  Our choice of notation is intentional in order to emphasize the connection between the Schr\"odinger problem and the multilinear operators discussed earlier.}
$$ -u_{\xi\xi}(\xi,x) + V(\xi) u(\xi,x) = x^2 u(\xi,x)$$
for some real number $x \neq 0$.  We are interested in the question of whether two linear independent solutions $u$ are both bounded for almost every $x$.
This would imply (among other things) that $[0,\infty)$ is an essential
support for the a.c. spectrum of the Schr\"odinger operator 
$-\partial_{\xi\xi} + V(\xi)$, a consequence that has recently been 
proved by completely different methods in \cite{deiftkillip}.
More quantitatively, we would like estimates on the maximal function 
$\sup_\xi |u(\xi,x)|$.

Formally, this eigenfunction equation has a solution
$$ u(\xi,x) = wkb(\xi) \sum_{n=0}^\infty (-1)^n T_{2n}(\check V, \ldots, \check V)(\xi,x)
+ \overline{wkb}(\xi) \sum_{n=1}^\infty (-1)^n T_{2n-1}(\check V, \ldots, \check V)(\xi,x)$$
where $T_n$ is the $n$-linear operator
$$ T_n(f_1, \ldots, f_n)(\xi,x) = (\frac{i}{2x})^n
\int_{\xi < \xi_1 < \xi_2 < \ldots < \xi_n}
\prod_{j=1}^n wkb(\xi_j)^{2(-1)^{n-j}} \hat f_j(\xi_j)\ d\xi_j$$
and the WKB phase is defined by
$$ wkb(\xi) := \exp(i x \xi - \frac{i}{2x} \int_0^\xi V).$$
See \cite{ck-2} for more details.
Thus to obtain bounds on this particular $u(\xi,x)$ and similarly on
all eigenfunctions $u$ it would suffice to obtain bounds on the operators 
$T_n$ which were decreasing sufficiently fast in $n$.  In the case when $V \in L^p$ for $p < 2$ this has been achieved in \cite{ck}; see also \cite{ck-2}.  However in the critical case $p=2$ it is not known whether the eigenfunctions are bounded for a.e. $x$.  (When $p>2$ boundedness can fail, see 
\cite{kls}).

As a model approximation let us replace the WKB phase and its various powers by the simpler phase $\exp(2\pi i x \xi)$.  Let us also only consider the limiting case $\xi = -\infty$ (instead of the supremum over all $\xi$).  
Then the operator $T_n$ simplifies to
$$ \tilde T_n(f_1, \ldots, f_n)(x) = (\frac{i}{2k})^n
\int_{\xi_1 < \xi_2 < \ldots < \xi_n}
\exp(2\pi i x(\xi_1 + \ldots + \xi_n)) \prod_{j=1}^n \hat f_j(\xi_j)\ d\xi_j.$$
When $n=1$ this operator is essentially the identity, while for $n=2$ this operator is essentially the bilinear Hilbert transform.  For $n=3$ the operator is essentially the trilinear operator $T$ in the introduction.  Thus in order to carry out the program of \cite{ck} in the endpoint case $p=2$ it is necessary\footnote{Of course, one must eventually re-instate the variable $\xi$ and then take suprema over $\xi$.  When $n=1$ this creates Carleson's maximal operator (which is of weak-type (2,2) \cite{carleson}), while for $n=2$ one obtains a hybrid of the Carleson operator and the bilinear Hilbert transform.  This operator will be considered in a later paper.}  (among other things) to bound $T$ on $L^2$.  This motivates the work of this paper and the sequel \cite{mtt:fourierbiest}.

\section{Appendix II: $T_{walsh,\P,\Q}$ and $T$}

In this rather informal section we briefly explain why $T_{walsh,\P,\Q}$ is the
natural Walsh model of $T_+$.  Here $T_+$ is defined like $T$ but with an
integration $0<\xi_1<\xi_2<\xi_3$ instead. This is a very minor modification.

For every positive dyadic interval $\omega$, let $\omega_l$ and $\omega_r$ denote the left and right halves of $w$ respectively.  The key observation is that for almost every $3$-tuple
$0<\xi_1<\xi_2<\xi_3$ there is a unique smallest positive dyadic interval
$\omega$ which contains all $\xi_j$ and either
$\xi_1,\xi_2\in\omega_{l}$ and $\xi_3\in\omega_{r}$, or
$\xi_1\in\omega_{l}$ and $\xi_2,\xi_3\in\omega_{r}$.  As a consequence we have the decomposition

\[T_+(f_1,f_2,f_3)=\]
\begin{equation}\label{unauna}
\sum_{\omega}
\int_{\xi_1<\xi_2<\xi_3;\,\xi_1,\xi_2\in\omega_{l};\,
\xi_3\in\omega_{r}}
\widehat{f}_1(\xi_1)\widehat{f}_2(\xi_2)\widehat{f}_3(\xi_3)
e^{2\pi ix(\xi_1+\xi_2+\xi_3)}\,d\xi_1 d\xi_2 d\xi_3
\end{equation}

\[+ \sum_{\omega}\int_{\xi_1<\xi_2<\xi_3;\,\xi_1\in\omega_{l};\,
\xi_2, \xi_3\in\omega_{r}}
\widehat{f}_1(\xi_1)\widehat{f}_2(\xi_2)\widehat{f}_3(\xi_3)
e^{2\pi ix(\xi_1+\xi_2+\xi_3)}\,d\xi_1 d\xi_2 d\xi_3.\]

We can rewrite the first term in \eqref{unauna} as
\[\sum_{\omega}\left(\int_{\xi_1<\xi_2;\,\xi_1,\xi_2\in\omega_{l}}
\widehat{f}_1(\xi_1)\widehat{f}_2(\xi_2)
e^{2\pi ix(\xi_1+\xi_2)}\,d\xi_1 d\xi_2 \right)
\left(\int_{\xi_3\in\omega_{r}}\widehat{f}_3(\xi_3)
e^{2\pi ix\xi_3}\,d\xi_3 \right)
\]


\[=\sum_\omega B(	f_1*\varphi_{\omega_l},
f_2*\varphi_{\omega_l}) (f_3*\varphi_{\omega_r})\]
where $\varphi_{\omega_l}$ and $\varphi_{\omega_r}$ are convolution kernels
adapted to the frequency intervals $\omega_l$ and $\omega_r$, and $B$
is the bilinear Hilbert transform.

Indeed, a Walsh model which reproduces (up to an inessential affine
transformation of the frequency support of the output) 
this time-frequency
behaviour is given by
\[
\sum_{Q\in\Q}\frac{1}{|I_Q|^{1/2}}
\langle B_{Q_1}(f_1,f_2),\phi_{Q_1}\rangle
\langle f_3,\phi_{Q_2}\rangle\phi_{Q_3} \]
as desired.
The second term in \ref{unauna} is discussed similarly.

\section{Appendix III: Proof of Proposition \ref{abstract}}\label{abstract-sec}

We now prove Proposition \ref{abstract}.  Fix the collection $\P$ and the collections $a^{(j)}_{P_j}$.  We adopt the shorthand
$$ SIZE_j := \size_j( (a^{(j)}_{P_j})_{P \in \P} ); \quad 
ENERGY_j := \energy_j( (a^{(j)}_{P_j})_{P \in \P} ).$$

We may of course assume that $a^{(j)}_{P_j}$ are always non-zero.  We begin by considering the contribution of a single tree:

\begin{lemma}[Tree estimate]\label{single-tree}
Let $T$ be a tree in $\P$, and $a^{(j)}_{P_j}$ be complex numbers for all $P \in T$ and $j=1,2,3$.  Then
$$
| \sum_{P \in T} \frac{1}{|I_P|^{1/2}} a^{(1)}_{P_1} a^{(2)}_{P_2} a^{(3)}_{P_3}|
\leq |I_T|
\prod_{j=1}^3 \size_j( (a^{(j)}_{P_j})_{P \in T} ).$$
\end{lemma}

\begin{proof}
Without loss of generality we may assume that $T$ is a 3-tree.  We then use H\"older to estimate the left-hand side by
$$ 
(\sum_{P \in T} |a^{(1)}_{P_1}|^2)^{1/2}
(\sum_{P \in T} |a^{(2)}_{P_2}|^2)^{1/2}
(\sup_{P \in T} \frac{|a^{(3)}_{P_3}|}{|I_P|^{1/2}}).$$
From Definition \ref{size-def} we have
$$ (\sum_{P \in T} |a^{(j)}_{P_j}|^2)^{1/2} \leq |I_T|^{1/2} 
\size_j( (a^{(j)}_{P_j})_{P \in T} )$$
for $j=1,2$.  Also, since the singleton tree $\{P\}$ is a $1$-tree with top $P$, we have
$$ \frac{|a^{(3)}_{P_3}|}{|I_P|^{1/2}} \leq
\size_3( (a^{(3)}_{P_3})_{P \in T} )$$
for all $P \in \T$.  The claim follows.
\end{proof}

To bootstrap this summation over $T$ to a summation over $\P$ we would like to partition $\P$ into trees $T$ for which one has control over $\sum_T |I_T|$.  This will be accomplished by 

\begin{proposition}\label{decomp}  Let $j = 1,2,3$, $\P'$ be a subset of $\P$, $n \in \Z$, and suppose that
$$ \size_j( (a^{(j)}_{P_j})_{P \in \P'} ) \leq 2^{-n} ENERGY_j.$$
Then we may decompose $\P' = \P'' \cup \P'''$ such that
\be{size-lower}
\size_j( (a^{(j)}_{P_j})_{P \in \P''} ) \leq 2^{-n-1}
ENERGY_j
\end{equation}
and that $\P'''$ can be written as the disjoint union of trees $\T$ such that
\be{tree-est}
\sum_{T \in \T} |I_T| \lesssim 2^{2n}.
\end{equation}
\end{proposition}

\begin{proof}
The idea is to initialize $\P''$ to equal $\P'$, and remove trees from $\P''$ one by one (placing them into $\P'''$) until \eqref{size-lower} is satisfied.

We assume by pigeonholing that we only have quartiles $P$ such that 
the length of $I_P$ is an even (odd) p[ower of $2$.

We describe the tree selection algorithm.  We shall need four collections $\T'_{i}, \T''_{i}$ of trees, where $i \neq j$; we initialize all four collections to be empty.  

Suppose that we can find an $i \neq j$ and a quartile $P^0 \in \P''$ such that
\be{plus}
\sum_{P \in \P'': P_i< P_i^0} |a^{(j)}_{P_j}|^2 \geq 2^{-2n-3} |I_{P^0}|.
\end{equation}
We may assume that $P^0_i$ is maximal with respect to this property and the tile order $<$.  Having assumed this maximality, we may then assume that $\xi_{P^0}$ is maximal if $i < j$, or minimal if $i > j$; here $\xi_{P^0}$ is the center of $\omega_{P^0}$.

We then place the $i$-tree
$$ \{ P \in \P'': P_i < P^0_i \}$$
with top $P^0$ into the collection $\T'_i$, and then remove all the quartiles in this tree from $\P''$.  We then place the $j$-tree
$$ \{ P \in \P'': P_j < P^0_j \}$$
with top $P^0$ into the collection $\T''_i$, and then remove all the quartiles in this tree from $\P''$.

We then repeat this procedure until there are no further quartiles $P^0 \in \P''$ which obey \eqref{plus}.

After completing this algorithm, none of the tiles $P^0$ in $\P''$ will obey  \eqref{plus}, so that \eqref{size-lower} holds for all $i$-trees in $\P''$.  (If the tree does not contain its top, we can break it up as the disjoint union of trees which do).  We then set 
$\T := \bigcup_{i \neq j} \T'_{i} \cup \T''_{i}$
and $\P' := \bigcup_{T \in \T} T$.  

It remains to prove \eqref{tree-est}.  Since the trees in $\T''_{i}$ have the same tops as those in $\T'_{i}$ it suffices to prove the estimate for $\T'_\pm$.  We shall only prove the claim for $i < j$, as the argument for $i > j$ is similar.

Fix $i < j$.  The key geometric observation is that the tiles
$$ \{ P_j:  P \in T \hbox{ for some } T \in \T'_i \}$$
are all pairwise disjoint.  Indeed, suppose that there existed $P \in T \in T'_i$ and $P' \in T' \in T'_i$ such that $P_j \neq P'_j$ and $P_j \cap P'_j 
\neq \emptyset$.  Without loss of generality we may assume that 
\be{G1}
I_P \supsetneq I_{P'}
\end{equation}
so that
$$\omega_{P_j} \subsetneq \omega_{P'_j}.$$  
From the nesting of dyadic intervals, and from the assumption that
two different scales differ at leats by a factor of $4$, this implies that 
$$\omega_{P_i} \subsetneq \omega_{P'_j}.$$
Since $\T'_i$ consists entirely of $i$-trees, we have 
$\omega_{P_{T,i}} \subset \omega_{P_i}$, thus
\be{G2}
\omega_{P_{T,i}} \subsetneq \omega_{P'_j}.
\end{equation}
On the other hand, since $T'$ is an $i$-tree, we have
$$\omega_{P_{T',i}} \subseteq \omega_{P'_i}.$$
Since $i < j$, we thus see that $\omega_{P_{T,i}}$ and $\omega_{P_{T',i}}$ are disjoint and that
$$ \xi_{P_{T,i}} > \xi_{P_{T',i}}.$$
Since we chose our trees $T$ in $\T'_i$ so that $\xi_{P_{T,i}}$ was maximized, this implies that $T$ was selected earlier than $T'$.  On the other hand,
from \eqref{G2} and the nesting of dyadic intervals we have
$$\omega_{P_{T,j}} \subsetneq \omega_{P'_j}$$
which implies from \eqref{G1} that
$$P'_j < P_{T,j}.$$
Thus $P'$ would have been selected for a tree in $\T''_i$ at the same time that $T$ was selected for $\T'_i$.   But this contradicts the fact that $P'$ is part of $T'$, and therefore selected at a later time for $\T'_i$.  This establishes the pairwise disjointness of the $P_j$.

From this disjointness and Definition \ref{size-def} we have
$$ \sum_{T \in \T'_i} \sum_{P \in T} |a^{(j)}_{P_j}|^2 \lesssim ENERGY_j.$$
From \eqref{plus} we therefore have
$$ \sum_{T \in \T'_i} 2^{-2n} |I_T| \lesssim ENERGY_j$$
as desired.
\end{proof}

From Proposition \ref{decomp} we easily have

\begin{corollary}\label{decomp-cor}  There exists a partition
$$ \P = \bigcup_{n \in \Z} \P_n$$
where for each $n \in \Z$ and $j = 1,2,3$ we have
$$ \size_j( (a^{(j)}_{P_j})_{P \in \P'} ) \leq \min(2^{-n} ENERGY_j, SIZE_j).$$
Also, we may cover $\P_n$ by a collection $\T_n$ of trees such that
$$ \sum_{T \in \T_n} |I_T| \lesssim 2^{2n}.$$
\end{corollary}

\begin{proof}
Since $\P$ is finite, we see that the hypotheses of Proposition \ref{decomp} hold for all $j=1,2,3$ if $n = -N_0$ for some sufficiently large $N_0$.  Set the $\P_n$ to be empty for all $n < -N_0$.  Now initialize $n = -N_0$ and $\P' = \P$.  For $j=1,2,3$ in turn, we apply Proposition \ref{decomp}, moving the quartiles in $\P'''$ from $\P'$ in $\P_n$ and keeping the tiles in $\P''$ inside $\P'$.  We then increment $n$ and repeat this process. Since we are assuming the $a^{(j)}_{P_j}$ are non-zero, every quartile must eventually be absorbed into one of the $\P_n$.  The properties are then easily verified.
\end{proof}

From Corollary \ref{decomp-cor} and Lemma \ref{single-tree} we see that
$$
| \sum_{P \in T} \frac{1}{|I_P|^{1/2}} a^{(1)}_{P_1} a^{(2)}_{P_2} a^{(3)}_{P_3}|
\leq |I_T|
\prod_{j=1}^3 \min(2^{-n} ENERGY_j, SIZE_j)$$
for all $T \in \T_n$.  Summing over all $T$ in $\T_n$ and then summing over all $n$, we obtain
$$
| \sum_{P \in T} \frac{1}{|I_P|^{1/2}} a^{(1)}_{P_1} a^{(2)}_{P_2} a^{(3)}_{P_3}|
\lesssim \sum_n 2^{2n}
\prod_{j=1}^3 \min(2^{-n} ENERGY_j, SIZE_j).$$
Without loss of generality we may assume that
$$ \frac{ENERGY_1}{SIZE_1} \leq \frac{ENERGY_2}{SIZE_2} \leq \frac{ENERGY_3}{SIZE_3}.$$
We may estimate the right-hand side as
\begin{align*}
 \sum_n \min(&2^n ENERGY_1 SIZE_2 SIZE_3, \\
& ENERGY_1 ENERGY_2 SIZE_3, \\
&2^{-n} ENERGY_1 ENERGY_2, ENERGY_3)
\end{align*}
which can be bounded by
$$ ENERGY_1 ENERGY_2 SIZE_3 \log(1 + \frac{ENERGY_3/SIZE_3}{ENERGY_2/SIZE_2}).$$
The claim then follows.

\end{document}